\documentclass[11pt]{article}
\usepackage{amsfonts}
\usepackage{mathrsfs}
\usepackage{graphics}
\usepackage{indentfirst}
\usepackage{amsthm}

\usepackage{cite}
\usepackage{latexsym}
\usepackage{amsmath}
\usepackage{amssymb}
\usepackage[dvips]{epsfig}
\usepackage{amscd}
\newcommand{\ba}{\begin{aligned}}
\newcommand{\ea}{\end{aligned}}
\newcommand{\be}{\begin{equation}}
\newcommand{\ee}{\end{equation}}
\newcommand{\bn}{\begin{eqnarray}}
\newcommand{\en}{\end{eqnarray}}
\newcommand{\bnn}{\begin{eqnarray*}}
\newcommand{\enn}{\end{eqnarray*}}
\newcommand{\ti}{\tilde}

\newcommand{\lm}{\lambda}
\renewcommand{\div}{ {\rm div }  }

\def\rr{\mathbb{R}^3}

\newcommand{\na}{\nabla }

\newcommand{\pa}{\partial}
\newcommand{\bi}{\bibitem}

\newcommand{\bl}{\begin{lemma}}
\newcommand{\el}{\end{lemma}}
\newcommand{\et}{\end{theorem}}

\newcommand{\te}{\theta}
\newcommand{\al}{\alpha}
\newcommand{\de}{\delta}
\newcommand{\ve}{\varepsilon}
\newcommand{\la}{\label}
\newcommand{\ka}{\kappa}
\newcommand{\om}{\Omega}

\newcommand{\ep}{\varepsilon}
\newcommand{\n}{\rho}

\renewcommand{\theequation}{\thesection.\arabic{equation}}
\def\la{\label}
\def\na{\nabla}

 \def\p{\partial}

\def\norm[#1]#2{\|#2\|_{#1}}
\def\nm[#1]#2{\|#2\|_{#1}}
\numberwithin{equation}{section}

 \usepackage{indentfirst}
\usepackage{graphics}
\usepackage{epsfig}

\newtheorem{claim}{\bf \t}[part]

\newtheorem{theorem}{Theorem}[section]

\newtheorem{lemma}[theorem]{Lemma}

\newtheorem{remark}{Remark}[section]
\newtheorem{definition}{Definition}[section]

\def\t{\theta}
\def\k{\kappa}

\def\m{\mu}

\def\l{\lambda}
\def\rr{\mathbb{R}^3}
\def\r{\rho}

\def\f{\frac}
\begin{document}

\title{Serrin-Type  Blowup Criterion for  Viscous, Compressible, and Heat Conducting Navier-Stokes and Magnetohydrodynamic  Flows
}

\author{
Xiangdi Huang,\thanks{ NCMIS, AMSS,   Chinese Academy of Sciences, Beijing
100190, People's Republic of China; Department of Pure and Applied Mathematics, Graduate School of Information Science and Technology,
Osaka University({\tt xdhuang@amss.ac.cn}). X. D. Huang is partially supported  by NNSFC 11101392.}
 \quad Jing Li \thanks{Institute of Applied Mathematics, AMSS, and Hua Loo-Keng Key Laboratory of Mathematics,
Chinese Academy of Sciences, Beijing 100190, People's  Republic of
China ({\tt ajingli@gmail.com}). J. Li is partially supported
 by the National Center for Mathematics and Interdisciplinary Sciences, CAS, and NNSFC 10971215 \& 11171326.}
 }
\date{ }
\maketitle

\begin{abstract}
This paper  establishes a blowup criterion for   the  three-dimensional  viscous, compressible, and heat conducting magnetohydrodynamic (MHD) flows.
 It is essentially shown that for the Cauchy problem and the initial-boundary-value one of   the three-dimensional compressible MHD flows with  initial density allowed to vanish, the strong or smooth solution exists globally if  the density is bounded from above and the velocity satisfies the Serrin's condition.  Therefore, if  the Serrin norm of the velocity  remains bounded, it is not possible for
other kinds of singularities (such as  vacuum states  vanish  or vacuum appears  in the non-vacuum region or even milder
singularities) to form before the density becomes unbounded.  This criterion is analogous to the well-known   Serrin's  blowup criterion for the  three-dimensional  incompressible Navier-Stokes equations, in particular, it is independent of the temperature and magnetic field  and is   just the same as that of the
barotropic compressible Navier-Stokes equations.
   As a direct application, it is shown that the same result  also holds for the strong or smooth solutions to the three-dimensional  full compressible Navier-Stokes system  describing the motion of a viscous, compressible, and heat conducting fluid.

\


\

AMS: 35Q35, 35B65, 76N10
\end{abstract}

\section{Introduction}

In this paper, we consider the system of partial differential
equations for the three-dimensional viscous,   compressible, and heat conducting magnetohydrodynamic (MHD) flows
in the Eulerian coordinates \cite{lif}
 \be  \label{a1}
\begin{cases}
\r_{t}+\mbox{div}(\r u)=0,\\
(\r u)_{t}+\mbox{div}(\r u\otimes
u)-\mu\Delta{u}-(\mu+\l)\nabla\mbox{div}u
+\nabla{P}= (\mbox{curl $H$})\times H,\\
c_v[(\r \t)_{t}+\mbox{div}(\r u\t
)]-\k\Delta{\t}+P\mbox{div}u=2\mu|\mathfrak{D}(u)|^2+\l(\mbox{div}u)^2 +
 {\nu}|\mbox{curl $H$}|^2,\\
H_t-\mbox{curl $(u\times H)$} = {\nu}\triangle H,\quad \mbox{div} H=0 ,
\end{cases}
\ee
where $t\ge 0$ is time, $x\in \Omega\subset \rr$ is the spatial coordinate, and
 $\n,  u=\left(u_1,u_2,u_3\right)^{\rm tr},  $ $\te, $ $P=R\r\t \,(R>0), $ and $ H=\left(H_1,H_2,H_3\right)^{\rm tr},$
represent respectively the fluid density,  velocity, absolute temperature, pressure, and magnetic field; $ \mathfrak{D}(u)$ is the
deformation tensor given by  \bnn \la{z1.2}   \mathfrak{D}(u)  =
\frac{1}{2}(\nabla u + (\nabla u)^{\rm tr}). \enn
 The constant
viscosity coefficients $\mu$ and $\lambda$  satisfy the physical
restrictions
 \be\la{a2} \mu>0,\quad 2 \mu + 3\lambda\ge 0.
\ee Positive constants $c_v,$  $\ka,$ and $\nu$ are respectively  the heat capacity,  the ratio of the heat conductivity coefficient over the heat capacity, and the magnetic diffusivity acting as a magnetic diffusion
coefficient of the magnetic field.

The equations (\ref{a1}) will be studied with  initial condition:
\be\la{a3}  (\rho,u,\te,H)(x,0)=(\rho_0, u_0,\te_0,H_0)(x),\quad x\in \om ,\ee and one of the  following  boundary conditions:

1) If $\om =\rr,$ for constant $\tilde \n\ge 0,$ $(\rho,u,\te,H)$ satisfies  the far field
condition:
\begin{equation}\label{bc1}
 (\r,u,H,\theta)(x,t)\rightarrow (\tilde \r,0,0, 0)
 ~~\mbox{as}~~|x|\rightarrow\infty;
\end{equation}

2) If $\om $ is a bounded smooth domain
in $\rr,$ $( u,\te,H)$ satisfies
   \be \label{bc2}u=0  ,\quad \frac{\pa \te}{\pa n}=0,\quad H=0 \quad\mbox{ on }\pa\om ,\ee  where $n=(n_1,n_2,n_3)$ is the unit outward normal to $\p\om .$

The compressible MHD system \eqref{a1}   is a combination of the compressible Navier-Stokes equations of fluid dynamics and Maxwell¡¯s equations
of electromagnetism. Indeed, the equations \eqref{a1}$_1$, \eqref{a1}$_2,$ and \eqref{a1}$_3$
 describe, respectively,   the conservation of mass, momentum, and energy.
In addition, it is well-known that the electromagnetic fields are governed by Maxwell¡¯s equations.
In magnetohydrodynamics, the displacement current can be neglected (\cite{lif}). As a
consequence, the equation  \eqref{a1}$_4$ is called the induction equation, and the electric field can be
written in terms of the magnetic field $H$ and the velocity $u,$
$$E=\nu\na\times H-u\times H.$$
Although the electric field $E$ does not appear in the compressible MHD system \eqref{a1}, it is indeed
induced according to the above relation by the moving conductive flow in the magnetic
field. In particular, when there is no electro-magnetic
effect, that is, $H\equiv 0,$  the  compressible MHD system \eqref{a1}  reduces to the following   full compressible Navier-Stokes system  describing the motion of a viscous, compressible, and heat conducting fluid:
\be  \label{n1}
\begin{cases}
\r_{t}+\mbox{div}(\r u)=0,\\
(\r u)_{t}+\mbox{div}(\r u\otimes
u)-\mu\Delta{u}-(\mu+\l)\nabla\mbox{div}u
+\nabla{P}=0,\\
c_v[(\r \t)_{t}+\mbox{div}(\r u\t
)]-\k\Delta{\t}+P\mbox{div}u=2\mu|\mathfrak{D}(u)|^2+\l(\mbox{div}u)^2.
\end{cases}
\ee

There is a considerable body of literature on the multi-dimensional full compressible Navier-Stokes system \eqref{n1} and
compressible MHD one \eqref{a1}  by physicists and mathematicians
because of their physical importance, complexity, rich phenomena, and mathematical
challenges; see \cite{vp,fjs,choe1,df,lif,huw1,feireisl1,Hof1,hlx4,L1,M1,Na,R,se1,X1} and the references cited therein. However, many physically important and mathematically fundamental
problems are still open due to the lack of smoothing mechanism and the strong nonlinearity. For example, although
  the local strong solutions to the compressible
MHD system \eqref{a1} with large initial data were respectively obtained by \cite{vp}
and \cite{fjs} in the cases that the initial density is strictly positive and that the  density is
allowed to vanish initially,
whether the unique local strong solution   can exist globally is an outstanding challenging open problem.

Therefore, it is important to study   the mechanism of blowup and
structure of possible singularities of strong (or smooth) solutions
to   the compressible MHD system \eqref{a1} and to the full compressible  Navier-Stokes one \eqref{n1}.
The pioneering work can be traced   to Serrin's criterion \cite{se2}  on the Leray-Hopf weak solutions to the three-dimensional incompressible Navier-Stokes
equations, which can be stated that if a weak solution $u$ satisfies
\begin{equation}\label{1.6} u\in L^s(0,T;L^r),\quad
\frac{2}{s}+\frac{3}{r}\leq 1,\quad 3<r\leq \infty,
\end{equation}
then it is regular. Later,
 He-Xin\cite{hex}   showed that the Serrin's criterion \eqref{1.6}  still holds even for the strong solution to the incompressible MHD equations.

 Recently, Huang-Li-Xin \cite{hlx}  extended  the Serrin's  criterion  \eqref{1.6}
to the barotropic compressible Navier-Stokes
equations  and  showed  that  if $T^*<\infty$ is the maximal
time of existence of a strong (or classical) solution $(\n,u)$, then    \be\la{ba1}\lim_{T\rightarrow T^*}\left(
\|{\text{div}u}\|_{L^1(0,T;L^{\infty})} + \| {
u}\|_{L^s(0,T;L^r)}\right) = \infty ,
 \ee  and \be\la{ba2} \lim_{T\rightarrow T^*}\left(
\|\rho\|_{L^\infty(0,T;L^{\infty})} + \| {
u}\|_{L^s(0,T;L^r)}\right) = \infty ,
 \ee    with $r$ and $s$ as in   (\ref{1.6}).  For more information  on the blowup criteria of barotropic compressible flow, we refer to
\cite{has,hlx,hlx1,H2,H4,wz} and the references therein. Later Xu-Zhang \cite{z} extended the results of \cite{hlx} to  the isentropic compressible MHD system and obtained that the same blow-up criterion \eqref{ba2} holds.

When it comes to the full compressible Navier-Stokes system \eqref{n1}, the problem is much more complicated.  Let $T^*<\infty$ be the maximal time of
existence of a strong (or classical) solution $(\n,u,\t)$ to the system \eqref{n1}. Besides \eqref{a2}, under the   condition that
\be\la{a6} 7\mu>\lambda, \ee
  Fan-Jiang-Ou \cite{J2}  obtained   that     \bnn
 \lim_{T\rightarrow T^*}(
\|{\theta}\|_{L^{\infty}(0,T;L^{\infty})} +
 \|{\nabla u}\|_{L^1(0,T;L^{\infty})})  = \infty .\enn
Recently,  under just the physical restrictions (\ref{a2}),    Huang-Li
\cite{hl} and Huang-Li-Xin \cite{hlx1}   established the
following  blowup criterion: \bnn \lim_{T\rightarrow
T^*}\left(
 \|{\theta}\|_{L^2(0,T;L^{\infty})}+\|{\mathfrak{D}( u)}
 \|_{L^1(0,T;L^{\infty})}\right)  = \infty ,\enn
where $ \mathfrak{D}(u)$ is the deformation tensor. Later, in the absence of vacuum,
Sun-Wang-Zhang \cite{wz1} showed that   \bnn \lim_{T\rightarrow
T^*}\left(\|{\theta}\|_{L^{\infty}(0,T;L^{\infty})} +
 \left\|{\left(\rho, \rho^{-1}\right)}\right\|_{L^{\infty}(0,T;L^{\infty})}\right)  = \infty,
\enn provided that  \eqref{a2} and  (\ref{a6}) both hold.
Very recently, under just the physical restrictions (\ref{a2}) and allowing the initial density to vanish,  Huang-Li-Wang \cite{HLW-1}  improved all the previous results \cite{wz1,hl,hlx1,J2} by obtaining that \eqref{ba1} still holds. It should be noted here that \eqref{ba2} is much stronger than \eqref{ba1} and that whether \eqref{ba2}  holds or not remains open.

For the compressible MHD system \eqref{a1}, let $T^*<\infty$ be the maximal time of
existence of a strong (or classical) solution $(\n,u,\t,H).$
Lu et al \cite{llya}  obtained   that     \bnn
 \lim_{T\rightarrow T^*}( \|{\rho}\|_{L^{\infty}(0,T;L^{\infty})}+
\|{\theta}\|_{L^{\infty}(0,T;L^{\infty})}+
 \|{\nabla u}\|_{L^4(0,T;L^2)})  = \infty ,\enn and  \bnn
 \lim_{T\rightarrow T^*}( \|{\div}u\|_{L^{\infty}(0,T;L^{\infty})}+
\|{\theta}\|_{L^{\infty}(0,T;L^{\infty})}+
 \|{\nabla u}\|_{L^4(0,T;L^2)})  = \infty,\enn while  Chen-Liu\cite{cl} showed that
  \bnn
 \lim_{T\rightarrow T^*}( \|{\nabla u}\|_{L^1(0,T;L^{\infty})}+
\|{\theta}\|_{L^{\infty}(0,T;L^{\infty})}
)  = \infty .\enn

The aim of this paper is  to  improve all the previous blowup
criterion  results on both the compressible MHD system \eqref{a1} and the full compressible Navier-Stokes one \eqref{n1} by allowing initial vacuum states, and by  describing the
blowup mechanism just in terms of  the Serrin-type criterion, \eqref{ba2}. Before stating our main result, we first  explain the
notations and conventions used throughout this paper. We denote
$$\int fdx=\int_{\om }fdx.$$
For $1\le p\le \infty$ and integer $k\ge 0$, the standard homogeneous and inhomogeneous Sobolev
spaces are denoted by:
   \bnn \begin{cases} L^p=L^p(\om ),\quad W^{k,p}=W^{k,p}(\om ),\quad
    D^{k,p}  =
    \left.\left\{u\in
L^1_{\rm loc}(\om )\,\right| {\nabla^k u}\in {L^p}  \right\},\\
  D^1_0   = \left. \left\{u\in L^6 \,\right|
  {\nabla u}\in{L^2} ,\,u = 0\mbox{ on }\pa \om  \right\},\quad H_0^1=L^2\cap D^1_0, \quad H^k=W^{k,2},\\ D^{2,2}_{0,n}   = \begin{cases}\left. \left\{\te\in D^{1,2}\cap D^{2,2} \,\right|
    {\na \te}\cdot{  n} = 0\mbox{ on }\pa \om  \right\},& \mbox{ for  bounded }\om     ,\\ D^1_0\cap D^{2,2},& \mbox{ for }\om =\rr. \end{cases}
 \end{cases}\enn
Then, the   strong solutions to the
  initial-boundary-value problem \eqref{a1}--\eqref{a3}    together with \eqref{bc1} or  \eqref{bc2} are defined as follows.
\begin{definition}[Strong Solutions] \label{def1.1}For $\ti\n\ge 0$ and $\ti\t=0,$
$(\r,u,\t,H)$ is called a strong solution to \eqref{a1} in
$\om \times (0,T)$, if for some $q_0>3$,
\begin{eqnarray}\nonumber
\begin{cases}
\r\geq 0,~~ \r-\tilde\r\in C([0,T];H^1\cap W^{1,q_0}),~~\r_t\in
C([0,T];L^2\cap L^{q_0}),\\
(H,u)\in C([0,T];D_0^1\cap D^{2,2})\cap L^2(0,T;D^{2,q_0}),\quad H\in C([0,T];H^2)\\ \t\ge 0,\quad
\t\in C([0,T]; D_{0,n}^{2,2})\cap L^2(0,T;D^{2,q_0}), \\
(H_t,u_t,\t_t)\in L^2(0,T;D^{1,2}), ~~
(H_t,\sqrt{\r}u_t,\sqrt{\r}\t_t)\in L^{\infty}(0,T; L^2),
\end{cases}
\end{eqnarray}
and $(\r,u,\t,H)$ satisfies both \eqref{a1} almost everywhere in $\om \times
(0,T)$ and  \eqref{a3} almost everywhere in $\om .$

\end{definition}

  Our main result can be stated as follows:

\begin{theorem}\label{thm1.1}
For $\tilde{q}\in(3,6]$, assume that the initial data $(\r_0\ge
0,u_0,\t_0\ge 0,H_0)$ satisfies
\begin{equation}\label{1.7}
\ba
 &\r_0-\tilde\r\in H^1\cap W^{1,\tilde{q}},\quad
  u_0 \in D_0^1\cap
D^{2,2},\quad \t_0 \in
D^{2,2}_{0,n}, \\&    \r_0|u_0|^4+ \r_0\t_0^2\in L^1,~~H_0\in H^1_0\cap  H^2,\quad\div  H_0  = 0,
\ea
\end{equation}
and the compatibility conditions
\be \label{1.8}
 -\m\Delta{u}_0-(\m+\l)\nabla {\rm div}u_0+R\nabla(\r_0\t_0)-( {\rm curl}H_0)\times H_0=\sqrt{\r_0}g_1,\ee
\be\label{1.9} \k\Delta\t_0+\f\m2|\nabla{u}_0+(\nabla{u}_0)^{tr}|^2+\l(\div u_0)^2+\nu|{\rm curl}  H_0  |^2=\sqrt{\r_0}g_2,
\ee
with $g_1,g_2\in L^2$. Let $(\r,u,\t,H)$ be the strong solution to the
  initial boundary value problem \eqref{a1}--\eqref{a3}    together with \eqref{bc1} or  \eqref{bc2}.
   If $T^\ast<\infty$ is the
maximal time of existence, then  for $r$ and $s$
as in \eqref{1.6},  \be \la{1.10}\lim_{T\rightarrow
T^*}(\norm[L^\infty(0,T;L^{\infty})]{\n} +
\norm[L^s(0,T;L^r)]{ u}) = \infty.
 \ee
\end{theorem}

If $ H\equiv  H_0\equiv 0,$ Theorem \ref{thm1.1} directly yields  the following Serrin-type blowup criterion for  the three-dimensional  full compressible Navier-Stokes system \eqref{n1}.

\begin{theorem}\label{thm1.2}
For constants $\tilde{q}\in(3,6]$ and $\ti\r\ge 0$, assume that  $(\r_0\ge
0,u_0,\t_0\ge 0)$ satisfies
\bnn
\ba
  \r_0-\tilde\r\in H^1\cap W^{1,\tilde{q}},\quad
  u_0 \in D_0^1\cap
D^{2,2},\quad \t_0 \in
D^{2,2}_{0,n}, \quad    \r_0|u_0|^4+ \r_0\t_0^2\in L^1,
\ea
\enn
and the compatibility conditions
$$ \label{n1.1}
 -\m\Delta{u}_0-(\m+\l)\nabla {\rm div}u_0+R\nabla(\r_0\t_0) =\sqrt{\r_0}g_1,$$ $$
 \label{n1.2} \k\Delta\t_0+\f\m2|\nabla{u}_0+(\nabla{u}_0)^{tr}|^2+\l(\div u_0)^2 =\sqrt{\r_0}g_2,
$$
with $g_1,g_2\in L^2$.
 Let $(\r,u,\t)$ be the strong solution to the
  full compressible Navier-Stokes system \eqref{n1}   together with
    \be \label{n1.3} (\n,u,\te)(x,0)=(\r_0,u_0,\te_0), \quad x\in \om ,\ee
 and either for $\om =\rr,$
    \be  \label{n1.4}  (\r,u,\te)\rightarrow (\ti\n,0,0)\mbox{ as }|x|\rightarrow \infty,  \ee or for a bounded smooth domain $ \om  \subset \rr,$\be \label{n1.5}  u=0,~~ \frac{\pa \te}{\pa n}=0 \mbox{ on }\pa\om   . \ee
     If $T^\ast<\infty$ is the
maximal time of existence, then  \be \la{fns2}\lim_{T\rightarrow
T^*}(\norm[L^\infty(0,T;L^{\infty})]{\n} +
\norm[L^s(0,T;L^r)]{ u}) = \infty,
 \ee with $r$ and $s$
as in \eqref{1.6}.
\end{theorem}

A few remarks are in order:

\begin{remark}
The conclusion  in  Theorem \ref{thm1.1} is somewhat surprising  since the criterion   \eqref{1.10} is independent of the temperature and magnetic fields and    just the same as those of barotropic compressible Navier-Stokes equations (\cite{hlx}).
\end{remark}

\begin{remark} In \cite[Theorem 1]{HLW-1}, we obtained that   \eqref{ba1} holds for the  Cauchy problem of  the full compressible Navier-Stokes system  \eqref{n1}. Thus,
  \be \la{l1w} \lim\limits_{T\rightarrow T^*} \|\div u\|_{L^1(0,T;L^\infty)} =\infty,\ee provided  that\be  \la{l4w}  \sup\limits_{0\le T\le T^*}\|u\|_{L^r(0,T;L^s)}<\infty ,   \ee     for $  r , s $
as in \eqref{1.6}.
It follows from   the continuity equation
$(\ref{n1})_1$ that for $t\in [0,T^*),$
\be  \la{l2w} \n(x,t)=\n_0(y(0;x,t))\exp\left\{ -\int_0^t\div u(y(s;x,t),s)ds\right\},\ee where
  $y(s;x,t)$ is
  the characteristic curve
defined by \bnn \frac{\rm d}{{\rm d}s}y=u(y,s),\quad y(t;x,t)=x  .\enn
   The combination of \eqref{l1w}  with \eqref{l2w} implies that there may hold for the density:

 1) The density remains bounded, that is, \be \la{l5w} \lim\limits_{T\rightarrow T^*} \|\r  \|_{L^\infty(0,T;L^\infty)}<\infty;\ee

 2) The density may concentrate, that is, \be  \la{l3w} \lim\limits_{T\rightarrow T^*} \|\r  \|_{L^\infty(0,T;L^\infty)} =\infty;\ee

3) Vacuum states may vanish:  There exists some $x_1\in \om $ and $x_1(t)$ satisfying  $\n_0(x_1)=0$ and $y(0;x_1(t),t)=x_1$ such that\be \la{l9w} \lim\limits_{t\rightarrow T^*} \r( x_1(t) ,t) \ge c_0>0;\ee

4) Vacuum states may appear   in the  non-vacuum region: There exists some $x_2\in \om $ and $x_2(t)$ satisfying  $\n_0(x_2)>0$ and $y(0;x_2(t),t)=x_2$ such that\be  \la{l6w} \lim\limits_{t\rightarrow T^*} \r( x_2(t) ,t) =0.\ee
 Then one may ask:  Which one or some  of \eqref{l5w}--\eqref{l6w} will happen?
 Theorem \ref{thm1.2} gives an answer to this question by  obtaining that  the density will concentrate  provided that  \eqref{l4w}   holds. In other words, if  the Serrin norm of the velocity  remains bounded, it is not possible for
other kinds of singularities (such as  vacuum states  vanish  or vacuum appears  in the non-vacuum region or even milder
singularities) to form before the density becomes unbounded. Moreover,  \eqref{fns2} still holds for the initial-boundary-value problem   \eqref{n1}  \eqref{n1.3}  \eqref{n1.5}.  Thus, Theorem \ref{thm1.2} greatly improves all the  previous blowup criterion for the full compressible Navier-Stokes system \eqref{n1} \cite{wz1,hl, H4,J2,HLW-1}.
\end{remark}

\begin{remark} If $\om $ is a bounded smooth domain of $\rr,$
Theorems \ref{thm1.1} and \ref{thm1.2} still hold when   the boundary condition $\na\te\cdot n|_{\pa\om }=0$ is  replaced  by   $\te|_{\pa\om }=0.$\end{remark}

\begin{remark}Theorems \ref{thm1.1} and \ref{thm1.2}  also hold respectively for  classical
 solutions to the     three-dimensional compressible MHD system \eqref{a1} and to the  full  compressible  Navier-Stokes one \eqref{n1}.
\end{remark}

We now comment on the analysis of this paper.

Let $(\rho,u,\t,H)$ be a strong
solution described in Theorem \ref{thm1.1}. Suppose that
\eqref{1.10} were false,  that is,
\begin{equation}\label{1.10-1-1}
\lim\limits_{T\rightarrow T^\ast}\left(
\|\rho\|_{L^\infty(0,T;L^\infty)}+ \|u\|_{L^s(0,T;L^r)}\right)
\leq M_0<+\infty.
\end{equation}
We want to show that
\bnn\la{result-1}
\ba
&\sup_{0\le t\le T^\ast}\left(\| \r-\ti\r\|_{H^1\cap
W^{1,\tilde{q}}}+\|\nabla{u}\|_{H^1}+\|\nabla\theta\|_{H^1}+\|H\|_{H^2}\right) \le C<+\infty.
\ea
\enn

Since  the  methods in all previous works \cite{llya,HLW-1,wz1,cl,J2,wz1} depend
crucially on either the $L^\infty_tL^\infty_x$-norm of the temperature $\theta$  or the $L^1_tL^\infty_x$-norm of the divergence of the velocity $\div u,$       some new ideas are needed to recover all the a priori estimates  just under   the assumption (\ref{1.10-1-1}) without any   a priori bounds on   the temperature, the magnetic field,  and  the divergence of the velocity. In fact,  we     prove (see Lemma \ref{lem3.2}) that a control of the Serrin norm
of the velocity and $L^\infty_{t}L^\infty_{x}$-norm of the density implies   a control on the  $L^\infty_t L^2_x$  norm of $\na u$. In order to obtain this control,  the key observation is that, instead of the temperature $\te,$ we treat  the total energy $E= c_v\t+\frac{1}{2}|u|^2,$  which in turn greatly reduces the difficulties arising from the high nonlinearities of  the temperature equation, \eqref{a1}$_3.$ Indeed,   multiplying   the equation of the conservative form of the total energy $E$ (see \eqref{en-1})   by $E$ yields that  to bound the $L^2_{ t}L^2_{x }$-norm of $\na E$ (see \eqref{zk3}), it is enough to control that of $|u||\nabla u|,$ which in fact can be reduced to the estimate of the $L^2_tL^6_x$-norm of $\nabla u$  (see \eqref{k0.9}). Then, to overcome  the difficulty  caused by the boundary when $\om $ is bounded, motivated by \cite{wz,has}, we decompose the velocity into two parts (see \eqref{lame} and \eqref{wl1}) which together with the $L^p$-estimate for the Lam\'e system  yield the desired bound on  the $L^2_tL^6_x$-norm of $\nabla u$ (see \eqref{3.14-3}). Finally, the   a priori estimates  on both the
$L^\infty_tL^p_x$-norm of the density gradient and  the $L^1_tL^\infty_x$-norm of the velocity
gradient   can be obtained   simultaneously  by
solving a logarithm Gronwall inequality based on a  logarithm  estimate for the Lam\'e system (see  Lemma \ref{lem2.3}) and the a priori estimates we
have just derived.

The rest of the paper is organized as follows: In the next section, we
collect some elementary facts and inequalities that will be needed
later. The main result, Theorem \ref{thm1.1},  is proved in Section \ref{sec3}.
\section{Preliminaries}
In this section, we recall some known facts and elementary
inequalities that will be used later.

First,  the following
existence and uniqueness of local strong solutions when the initial
density may not be positive and may vanish in an open set  can be
proved in a similar way as in \cite{choe1} (cf. \cite{fjs}).
\begin{lemma}\label{lem2.1}
Assume that the initial data $(\r_0\ge 0,u_0,\t_0\ge 0,H_0)$ satisfy
\eqref{1.7}-\eqref{1.9}. Then there exists a positive time
$T_1\in(0,\infty)$ and a unique strong solution $(\r,u,\t,H)$ to the
  initial-boundary-value problem \eqref{a1}--\eqref{a3}    together with \eqref{bc1} or  \eqref{bc2} on $\om \times
(0,T_1]$.
\end{lemma}

Next, the following well-known Sobolev inequality will be used later
frequently (see \cite{nir}).
\begin{lemma}\label{lem2.2}
For $p\in(1,\infty)$ and $q\in(3,\infty)$, there exists a generic
constant $C>0$, which depends only on $p,~q$ such that for $f\in D_0^1$
and $g\in L^p\cap D^{1,q}$, we have
\begin{equation}\label{2.1}
\|f\|_{L^6}\leq C\|\nabla f\|_{L^2},~~\|g\|_{L^\infty}\leq
C\|g\|_{L^p}+C\|\nabla g\|_{L^q}.
\end{equation}
\end{lemma}

Finally, we consider  the following Lam\'e   system \be\la{u89} -\mu \Delta v(x)-(\mu+\lm)\na \div v(x)=  f(x) ,\quad   x\in\om , \ee where $v=(v_1,v_2,v_3),$   $f=(f_1,f_2,f_3),$ and $\mu,\lm$ satisfy \eqref{a2}. The system \eqref{u89} is imposed on one of the following boundary conditions:

 1) Cauchy problem:  $\om =\rr,$   and \be\la{b1} v(x)\rightarrow 0,\,\,\mbox{ as } |x|\rightarrow \infty;\ee

2) Dirichlet problem: $\om $ is a bounded smooth domain
in $\rr,$ and \be \la{b2} v=0 \mbox{ on }{\p\om }.\ee

The following  logarithm  estimate for the Lam\'e system \eqref{u89} will be used   to estimate $\|\nabla
u\|_{L^\infty}$ and $\|\nabla\r\|_{L^2\cap L^q}.$

\begin{lemma}\label{lem2.3} Let $\mu,\lm$ satisfy \eqref{a2}.
Assume that  $f=\div g$ where $g=(g_{kj})_{3\times 3} $ with $g_{kj}\in L^2\cap L^r\cap D^{1,q}$ for    $k,j=1,\cdots,3, $  $r\in (1,\infty), $  and $q\in (3,\infty).$  Then the  Lam\'e system  \eqref{u89} together with \eqref{b1} or  \eqref{b2}   has a unique solution $v\in   D_0^1\cap D^{1,r}\cap D^{2,q},$ and  there exists a generic positive constant $C$ depending only on $\mu,\lm,q, $ and $r$ (besides $\om $ when $\om $ is bounded) such that
\be \la{b6}\|\na  v\|_{L^r}\le C\|g\|_{L^r},\ee
and \be\la{ww7} \|\na v\|_{L^\infty}\le C \left(1+\ln
(e+\|\na g\|_{L^q})\|g\|_{L^\infty}+ \|g\|_{L^r}\right).\ee

\end{lemma}

{\it Proof.} First,  if $\om =\rr,$  direct calculations show that $v=(v_1,v_2,v_3)$ with $$v_j=\frac{1}{2\mu+\lm}(-\Delta)^{-1}\pa_k g_{kj},\quad j=1,\cdots, 3$$ is the unique solution to the Cauchy problem \eqref{u89} \eqref{b1} and satisfies \eqref{b6}.

Then, if $\om $ is a bounded smooth domain of $\rr,$ it follows from  \cite{sol2} that the Dirichlet problem \eqref{u89}   \eqref{b2}  is of Petrovsky type. In Petrovsky's systems, roughly speaking, different equations
and unknowns have the same ``differentiability order", see \cite[p.126]{sol1}.  We also recall that Petrovsky's systems are an important subclass of Agmon-Douglis-Nirenberg (ADN) elliptic systems(\cite{adn1}), having the same good properties of self-adjoint ADN
systems. It follows from
Solonnikov \cite[Theorem 1.1]{sol1} and \cite[Theorem 5.1]{sol2} that  the solution $v$ to  the system \eqref{u89} together with \eqref{b2}  can be represented as
  \be\la{b7} v_i(x)=\int G_{ij}(x,y)f_j(y)dy,\quad\mbox{ for all }x\in \om ,\ee by means of the Green function
 $G_{ij}=G_{ij}(x,y)\in C^\infty(\om \times\om \setminus D)$
with $D\equiv\{(x,y)\in \om \times\om |x=y\}$ which satisfies that for every multi-indexes $\al=(\al_1,\al_2,\al_3)$ and $\beta=(\beta_1,\beta_2,\beta_3)$ there is a
constant $C_{\al,\beta}$ such that for all $(x,y)\in \om \times\om \setminus D,$ and $i,j=1,\cdots, 3, $  \be \la{ww1}|\pa_x^\al\pa_y^\beta G_{ij}(x,y)|\le C_{\al,\beta}|x-y|^{-1-|\al|-|\beta|},\ee
where $|\al|=\al_1+\al_2+\al_3$ and $\beta=\beta_1+\beta_2+\beta_3.$
Moreover, the estimate \eqref{b6} is standard.

Finally, it remains to prove \eqref{ww7}. We will only deal with the Dirichlet problem \eqref{u89}   \eqref{b2}, since the same procedure holds for the Cauchy problem \eqref{u89}   \eqref{b1}. Motivated by    Beale-Kato-Majda \cite{bkm}, we introduce a small parameter  $\de\in(0,1]$ which depends on $v,$ and which will be fixed later. Using $\de,$ we define a cut-off function $\eta_\de(s)$
satisfying $\eta_\de(s)=1$ for $0\le s<\de,\eta_\de(s)=0$ for
$s>2\de,$ and $| \eta_\de^{(k)}(s)|\le C\de^{-k}.$
 It thus follows from \eqref{b7} that
\bnn\ba v_i(x)&=\int\left( \eta_\de (|x-y|)+ (1-\eta_\de (|x-y|))\right)G_{ij}(x,y)\pa_k  g_{kj}(y)dy\\ &=\int \eta_\de (|x-y|)G_{ij}(x,y)\pa_k  g_{kj}(y)dy+\int  \pa_{y_k} \eta_\de (|x-y|) G_{ij}(x,y) g_{kj}(y)dy\\&\quad-\int (1-\eta_\de (|x-y|))\pa_{y_k}G_{ij}(x,y)  g_{kj}(y)dy,\ea\enn
where in the second equality we have used integrations by parts due to the fact that $G_{ij}(x,y) |_{\pa\om }=0$ for each $x\in \om .$ Hence, we have
\be\la{ww2}\ba  |\na v(x)| &\le C\int\left(|\eta_\de'||G_{ij} |+\eta_\de|\na_xG_{ij}|\right)|\na g|dy\\&\quad +C\int\left(|\eta_\de''||G_{ij} | +|\eta_\de'||\na_xG_{ij} |+|\eta_\de'||\na_yG_{ij} |\right)|g|dy
\\&\quad+C\int(1-\eta_\de ) |\na_x\na_yG_{ij} ||g|dy.\ea\ee   Each term on the
right-hand side of (\ref{ww2}) can be estimated by (\ref{ww1}) as
follows: \be\ba & \int\left(|\eta_\de'||G_{ij} |+\eta_\de|\na_xG_{ij} |\right)|\na g|dy\\ &\le  C\int_{\om \cap \{y||x-y|<2\de\}}\left(\de^{-1}|x-y |^{-1}+ |x-y |^{-2}\right)|\na g|dy\\&\le C\left(\de^{-q/(q-1)}\int_0^{2\de} s^{-q/(q-1)}s^2ds+\int_0^{2\de} s^{-2q/(q-1)}s^2ds\right)^{(q-1)/q} \|\na g\|_{L^q}\\&\le C\de^{(q-3)/q}\|\na g\|_{L^q},\ea\ee
\be\ba &\int\left(|\eta_\de''||G_{ij} | +|\eta_\de'||\na_xG_{ij} |+|\eta_\de'||\na_yG_{ij} |\right)|g|dy\\ &\le C\int_\de^{2\de}\left(\de^{-2} s^{-1} +\de^{-1}s^{-2}\right)s^2ds\|g\|_{L^\infty}\\ &\le C \|g\|_{L^\infty},\ea\ee
and
\be\la{ww5}\ba &\int(1-\eta_\de ) |\na_x\na_yG_{ij} ||g|dy\\&\le C\left(
\int_{\om \cap\{y|\de\le |x-y|\le 1\}} +\int_{\om \cap\{y|  |x-y|>1\}}\right) |x-y|^{-3}|g(y) |dy \\&\le C \int_\de^1 s^{-3}s^2ds\|g\|_{L^\infty} +C\left(\int_1^\infty s^{-3r/(r-1)}s^2ds\right)^{(r-1)/r}
\|g\|_{L^r}\\&\le - C \|g\|_{L^\infty} \ln \de+C
\|g \|_{L^r}.\ea\ee
 It follows from (\ref{ww2})-(\ref{ww5})
that
\be\la{ww6} \|\na v\|_{L^\infty}\le
C\left(\de^{(q-3)/q}\|\na g\|_{L^q}+(1-\ln \de)\|g\|_{L^\infty}+ \|g\|_{L^r}\right).\ee Set
$\de=\min\left\{1,\|\na g\|_{L^q}^{-q/(q-3)}\right\}.$ Then
(\ref{ww6}) becomes (\ref{ww7}). We finish the proof of Lemma \ref{lem2.3}.

\section{\la{sec3}Proof  of Theorem \ref{thm1.1} }

 Before proving Theorem \ref{thm1.1}, we  state some a priori estimates under the condition (\ref{1.10-1-1}). First, we have
\begin{lemma}\label{lem3.1-1}
Under the condition \eqref{1.10-1-1}, it holds that for  $q\in[2,12]$ and $0\leq
T<T^\ast$,
\begin{equation}\label{3.3-1}
  \norm[L^\infty(0,T;L^q)]{H} + \int_0^T\int|H|^{q-2}|\nabla H|^2dx\le C,
\end{equation}
where (and in what follows)
$C $ and $C_i(i=1,\cdots, 6) $     denote  generic constants depending only on $
  M_0,$ $ \mu,$ $\lambda, $  $R,\kappa,$ $c_v,$  $T^*,\nu,$  and the
initial data (besides $\om $ for bounded $\om $).
\end{lemma}

 {\it Proof.} Similar to \cite{z}, multiplying $(1.1)_4$ by $q|H|^{q-2}H$ and integrating the resulting equation over $\om $ yield that
\be\la{hh1h}
\ba
& \f{d}{dt}\int|H|^qdx + \nu\int\left(q|H|^{q-2}|\nabla H|^2 +
q(q-2)|H|^{q-2}|\nabla|H||^2\right)dx\\
& = -\int q|H|^{q-2}\left(H\cdot\nabla H\cdot u - \f{q-1}{2}u\cdot\nabla|H|^2\right)dx\\
&\quad - \f{q(q-2)}{2}\int|H|^{q-4}(H\cdot\nabla|H|^2)(u\cdot H)dx\\
& \le \f{\nu}{2 }\int q|H|^{q-2}|\nabla H|^2 dx + Cq^2\int|u|^2|H|^qdx\\
& \le \f{\nu}{2 }\int q|H|^{q-2}|\nabla H|^2 dx + C \|u\|_{L^r}^2\||H|^{q/2}\|_{L^2}^{2(r-3)/r}\||H|^{q/2}\|_{L^6}^{6/r}
 \\ &\le \f{\nu}{2 }\int q|H|^{q-2}|\nabla H|^2 dx+C\de\norm[L^2]{\nabla |H|^{q/2}
}^2 + C(\de)(1+\norm[L^r]{u}^s)\norm[L^q]{ H }^q.\ea
\ee
   Choosing $\de$ suitably
small in \eqref{hh1h}, we obtain  \eqref{3.3-1} directly after using Gronwall's inequality and   $(\ref{1.10-1-1}).$  We thus finish the proof of Lemma  \ref{lem3.1-1}.

Then, we derive the following key  estimate on the  specific  energy  $E$ defined by \be\la{energy}
E\triangleq  c_v\te+\frac{|u|^2}{2}.\ee

\begin{lemma}\label{En}
  Under the condition \eqref{1.10-1-1}, it holds that
\be\la{zk3}
\ba
  c_v\frac{d}{dt}\int\r E^2 dx+ {\k} \norm[L^2]{\nabla E}^2
 \le  &C  \int|u|^2\left(\n E^2   +  |\na u|^2\right)dx+C_1  \|\na H \|_{H^1}^2 \\ & +C \|\na  u \|_{L^2}^2+C  \int\n E^2 dx .
\ea
\ee
\end{lemma}

{\it Proof.}
First,
  it follows from $\eqref{a1} $   that $E$ satisfies
\begin{equation}\label{en-1}(\r E)_t+\mbox{div}(\r Eu) -\frac{\k}{ c_v}   \Delta E =\div F -H^iH^j\pa_iu^j+\frac12|H|^2\div u+\nu|\mbox{curl} H|^2, \ee
with
$$
F  \triangleq  \frac{\mu-\k  c_v^{-1}}{2}\na(|u|^2) +\m u\cdot\nabla u +\l u\mbox{div}u    - Pu +( u\cdot H)H -\f12  |H|^2u .
$$
Next, applying standard maximum principle to $\eqref{a1}_3 $ together with
   $\theta_0\ge 0$  (c.f. \cite{J2,feireisl1}) shows  $$\label{3.07} \inf_{\mathbb{R}^3\times [0,T]}\theta(x,t)\ge 0.$$
      Multiplying \eqref{en-1} by $ c_vE$ and integrating the resulting equality over $\om ,$ we obtain after integration by parts and using $\eqref{a1}_1$ that \be\la{k0.4}\ba &\frac{ c_v}{2} \frac{d}{dt}\int \n E^2dx+ \ka \int |\na E|^2dx \\& \le C\int \left(|u| |\na u| +  \n\te| u| \right)|\na E| dx+C\int \left(|u||H|^2|\na E| +  |\na u||H|^2  E\right) dx\\&\quad+C \int E|{\rm curl} H|^2dx.\ea\ee

We estimate each term on the right-hand side of \eqref{k0.4} as follows:

First, Holder's inequality gives
\be\la{k0.3}\ba \int \left(|u| |\na u| +  \n\te| u| \right)|\na E| dx\le \eta\|\na E\|_{L^2}^2+C(\eta)\int\left( |u|^2|\na u|^2+ \n E^2|  u|^2\right) dx.\ea\ee

Next, if $\om =\rr,$ Sobolev's inequality gives that there exists a universal constant $C$ such that
\be\la{k0.1} \|E\|_{L^6}\le C\|\na E\|_{L^2}.\ee  If $\om $ is a bounded  smooth domain in $\rr,$
the Poincar\'e-type inequality     (\cite[Lemma 3.2]{feireisl1}) shows there exists a generic positive constant $C$ which also depends on $\om $ such that
\bnn \|E\|_{L^6}\le C\|\n^{1/2}E\|_{L^2}+C\|\na E\|_{L^2},\enn which combined with \eqref{k0.1} implies
\be \la{k0.2} \|E\|_{L^6}\le C\|\n^{1/2}E\|_{L^2}+C\|\na E\|_{L^2}.\ee It thus follows from Holder's inequality, \eqref{k0.2}, \eqref{3.3-1},  and \eqref{2.1} that
\be\la{k0.5}\ba &\int \left(|u||H|^2|\na E| +  |\na u||H|^2  E\right) dx\\&\le \|u\|_{L^6}\|H\|_{L^6}^2\|\na E\|_{L^2} +\|\na u\|_{L^{2}}\|H\|_{L^6}^2\|  E\|_{L^6}  \\ & \le \eta \|\na E\|_{L^2}^2+C(\eta)\|\n^{1/2}E\|_{L^2}^2+C(\eta)\|\na u\|_{L^2}^2 . \ea\ee

Finally, integration by parts together with  \eqref{k0.2}  yields
\be\la{k0.6}\ba \int E|{\rm curl} H|^2dx&\le C\int |\na E||\na H||H| dx+ C\int |  E||\na^2 H||H| dx\\ &\le C\|\na E\|_{L^2}\|\na H\|_{L^6}\|H\|_{L^3}+ C\|  E\|_{L^6}\|\na^2 H\|_{L^2}\|H\|_{L^3} \\ & \le \eta \|\na E\|_{L^2}^2+C(\eta)\|\n^{1/2}E\|_{L^2}^2+C(\eta)\|\na H\|_{H^1}^2 .\ea\ee

Putting \eqref{k0.3}, \eqref{k0.5}, and \eqref{k0.6} into  \eqref{k0.4},
 we obtain \eqref{zk3} after choosing $\eta$ suitably small.
 The proof of Lemma $\ref{En}$ is completed.

Then, we derive the following crucial  estimate on the $L^\infty(0,T;
L^2)$-norm of $\nabla u.$
\begin{lemma}\label{lem3.2}
Under the condition \eqref{1.10-1-1}, it holds that for $0\leq
T<T^\ast$,
\be\label{3.4}\ba &\sup_{0\leq t\leq T}
\int\left((\r-\ti\r)^2+\r \t^2 +|\nabla u|^2+|\nabla H|^2\right)
dx\\&+\int_{0}^{T}\int\left(|\nabla \t|^2+\r|\dot{u}|^2+|H_t|^2+|\na^2H|^2\right)dxdt \leq C.\ea
\ee
\end{lemma}

{\it Proof.}
 First, multiplying $(\ref{a1})_2$ by $u_t$ and integrating the resulting
equation over $\om $ show that
\be \label{3.14}\ba
& \f12\f{d}{dt}\int\left(\m|\nabla u|^2+(\m+\l)(\mbox{div}u)^2\right)dx+\int\r|\dot{u}|^2dx\\
&= \int\r\dot{u}\cdot(u\cdot\nabla){u} dx+\int P\mbox{div}u_tdx
-\f12\int (\na |H|^2-2\div(H\otimes H))\cdot u_tdx\\
&  \le\f14\int\r|\dot{u}|^2dx+C\int\n|u|^2|\nabla u|^2dx+\f{d}{dt}\int P\mbox{div}udx
\\& \quad -\int
P_t\mbox{div}udx+\frac{1}{2}\int(|H|^2\div u_t-2H\cdot\nabla u_t\cdot H)dx .
\ea\ee

Then, we will estimate the last two terms on the right-hand side of \eqref{3.14}.

On the one hand, to overcome  the difficulty  caused by the boundary,  motivated by \cite{wz,has}, we decompose the velocity into two parts. It follows from Lemma \ref{lem2.3} that for any $t\in [0,T], $ there exists a unique  $v(t,\cdot) \in D^1_0\cap D^{2,2}\cap D^{2,\ti q}$ satisfying
\be\label{lame}
 \mu\triangle v  + (\mu+\lambda)\nabla\mbox{div}v  = \nabla P  ,
\ee which together with \eqref{b6}  yields that \be \label{3.n3}\|\na v \|_{L^p}\le C\|P \|_{L^p}\le C\|\n E \|_{L^p}, \mbox{ for } p\in [2,6],\,\, t\in [0,T],\ee and that
\be \label{3.n1} \ba
-\int P_t {\rm div} v dx
& = -  \int (\mu\na v_t\cdot\na v + (\mu+\lambda) \mbox{div}v_t\div v) dx \\
& =- \frac12 \frac{d}{dt} \int \left(\mu | \na  v |^2 +(\mu+\lambda)(\div v)^2\right)dx.
\ea \ee
Denoting by \be\la{dw} w\triangleq u-v,\ee   we have  $w \in D^1_0\cap D^{2,2}\cap D^{2,\ti q},$  for a.e. $t\in [0,T].$ Moreover, for a.e. $t\in [0,T],$ $w$  satisfies \be\la{wl1}
 \mu\triangle w + (\mu+\lambda)\nabla\mbox{div}w = \n \dot u+H\times ({\rm curl}H),
\ee
which together with the standard $L^2$-estimate for elliptic system  gives \be \label{3.n2}\ba \|\na w\|_{L^6}+\|\na^2 w\|_{L^2}&\le C \|\n \dot u\|_{L^2}+C\||H||\na H| \|_{L^2} \\&\le C \|\n \dot u\|_{L^2}+C\|H\|_{L^6}\|\na H  \|_{L^2}^{1/2}\|\na H  \|_{L^6}^{1/2} \\&\le C \|\n \dot u\|_{L^2}+ C\|\na H  \|_{L^2}^{1/2}\|\na H  \|_{H^1}^{1/2} ,\ea\ee due to \eqref{3.3-1}.
 It follows from \eqref{energy} and \eqref{en-1} that
\be \label{3.m1} \ba
   -\int P_t  {\rm div} w dx    & =-\frac{R}{c_v}  \int (\rho E)_t {\rm div} w dx+ \frac{R}{2c_v}  \int (\rho |u|^2)_t {\rm div} w dx   \\
 & \le C   \int \left( \rho E |u| + |\na E|+|u||\na u|+|u||H|^2\right) |\na^2w| dx\\&\quad  +C\int  |H|^2|\na u||\na w|dx  -\frac{R\nu}{c_v}\int |\mbox{curl} H|^2  \div w dx
\\&\quad   - \frac{R}{2c_v}\int\left( {\rm div} (\rho u) |u|^2 {\rm div} w         - 2 \rho u\cdot u_t   {\rm div} w \right)dx   =\sum\limits_{i=1}^4I_i.
\ea \ee

Cauchy's and Sobolev's inequalities together with \eqref{3.3-1} yield that
\be\ba I_1+I_2\le &\eta\left(\|\na^2w\|_{L^2}^2+\|\na w\|_{L^6}^2\right)\\&+C(\eta)\int\left(\n^2E^2|u|^2+|\na E|^2+|u|^2|\na u|^2+|\na u|^2\right)dx.\ea\ee

 Similar to \eqref{k0.6}, integration by parts leads to
 \be \ba I_3\le & C\int\left(|\na H||H||\na^2w|+|\na^2H||H||\na w|
  \right)dx\\ \le &C\left(\|\na H\|_{L^6}\|\na^2w\|_{L^2}+\|\na^2H\|_{L^2}\|\na w\|_{L^6}\right)\|H\|_{L^3}\\ \le &\eta\left(\|\na^2w\|_{L^2}^2+\|\na w\|_{L^6}^2\right)+C(\eta)\|\na H\|_{H^1}^2,\ea \ee
where in the last inequality we have used \eqref{3.3-1}.

 Integration by parts also gives
\be\ba I_4&\le C\int \n |u|^3|\na^2w|dx+C\int \left(\n|u|^2|\na u|+\n |u||\dot u|\right)|\na w|dx\\ & \le  C(\eta)\int \left(\n E^2|u|^2+\n|u|^2|\na u|^2+\n |u|^2|\na v|^2\right) dx\\&\quad +\eta\|\na^2w\|_{L^2}^2+\eta\int\n |\dot u|^2dx.\ea\ee

On the other hand, direct calculations show
\be\la{zk5}
\ba
  &  \int( |H|^2\div u_t-2H\cdot\nabla u_t\cdot H)dx\\
& = \frac{d}{dt}\int( |H|^2\div u-2H\cdot\nabla u\cdot H)dx\\
&\quad -2 \int(H\cdot H_t\div u- H_t\cdot\nabla u\cdot H- H\cdot\nabla u\cdot H_t)dx\\
& \le \frac{d}{dt}\int( |H|^2\div u-2H\cdot\nabla u\cdot H)dx+
C\|{H_t}\|_{L^2}^2+C\||H||\na u|\|_{L^2}^2.
\ea
\ee
Substituting \eqref{3.n1}   and  \eqref{3.m1}--\eqref{zk5}  into \eqref{3.14}, we obtain after using \eqref{3.n2}  and choosing $\eta$ suitably small that
\be\la{zk1}\ba &\frac{d}{dt}\int\Phi dx +\int\n|\dot u|^2dx\\&\le  C \int |u|^2\left(\n E^2 +|\na u|^2  +   |\na v|^2\right)dx+ C \|\nabla u\|_{L^2}^2 \\&\quad +
C_2\left( \|\na  E \|_{L^2}^2+\|{H_t}\|_{L^2}^2+ \|\na H \|_{H^1}^2 + \||H||\na u|\|_{L^2}^2 \right),\ea\ee
where  \bnn\ba\Phi\triangleq & \m|\nabla u|^2+(\m+\l)(\mbox{div}u)^2-2P\div u+\m|\nabla v|^2+(\m+\l)(\mbox{div}v)^2\\&
 -|H|^2\div u+2H\cdot\nabla u\cdot H \ea\enn satisfies \be\la{p23} \Phi\ge \frac{\mu}{2} |\na u|^2 -C_3\n E^2-C|H|^4.\ee

Next, it follows from \eqref{a1} that for $r,s$ as in \eqref{1.6}
\be\la{a3.14-4}
\ba
&\nu \frac{d}{dt}\|\na  H\|_{L^2}^2  +\|H_t\|_{L^2}^2 +\nu^2 \|\Delta H\|_{L^2}^2  \\&=\int | H_t-\nu \Delta H|^2dx \\
& \le  C( \norm[L^2]{|H||\nabla u|}^2 + \norm[L^2]{|u||\nabla H|}^2)\\&\le C\|H\|_{L^6}^2\|\na u\|_{L^2}\|\na u\|_{L^6}+C  \norm[L^r]{u}^2 \norm[L^2]{\nabla H}^{2(r-3)/r} \norm[L^6]{\nabla  H}^{6/r}\\&\le \eta\|\na u\|_{L^6}^2+C(\eta) \|\na u\|_{L^2}^2 + C(\ve)   (1+\norm[L^r]{u}^{ {s} })  \norm[L^2]{\nabla  H}^2  +\ve\norm[L^2]{\na^2 H}^2 .
\ea
\ee
Noticing that the standard $L^2$-estimate of elliptic system  gives
 \bnn\ba
 \|\nabla^2H\|_{L^2}\le C_4  \|\Delta H\|_{L^2} ,
\ea\enn
after choosing $\ve$ suitably small, we deduce from \eqref{a3.14-4}  that for any $\eta\in(0,1),$
\be\la{zk2}
\ba
&  4\nu\frac{d}{dt}\|\nabla  H\|_{L^2}^2  +4 \|H_t\|_{L^2}^2 + 2\nu^2C_4^{-1}\left(\|\na  H\|_{H^1}^2 +\norm[L^2]{|H||\nabla u|}^2\right) \\&\le C\eta\|\na u\|_{L^6}^2+C(\eta) \|\na u\|_{L^2}^2 + C   (1+\norm[L^r]{u}^{ {s} })  \norm[L^2]{\nabla  H}^2  .
\ea
\ee

Then,  adding \eqref{zk3} multiplied by $C_5\triangleq C_3c_v^{-1}+(C_2+2)\ka^{-1}  $ and  \eqref{zk2}  by $C_6\triangleq  (1 +C_4\nu^{-2})(C_2+C_1C_5+2) $ to \eqref{zk1}, we obtain  that
\be\la{p329}\ba \frac{d}{dt}& \int\left( \Phi +C_5\n E^2 +4 C_6\nu|\na H|^2\right)dx \\& + \|\na E\|_{L^2}^2+ \frac12\int\n|\dot u|^2dx+ \|H_t\|_{L^2}^2+\|\na H\|_{H^1}^2\\ \le &C  \int|u|^2\left(\n E^2+|\na u|^2 +|\na v|^2 \right)dx +C  (1+\norm[L^r]{u}^{ {s} }) \norm[L^2]{\nabla H}^2  \\& +C\eta\|\na u\|_{L^6}^2 +C (\eta)\int \left(\n E^2  +  |\na u|^2 \right)dx   .\ea\ee
Holder's inequality together with \eqref{3.n3} yields that
\be\la{k0.9}\ba &\int|u|^2\left(\n E^2+|\na u|^2+|\na v|^2  \right)dx \\   & \le C\norm[L^r]{u}^2  \left(\|\n^{1/2}E\|_{L^{2r/(r-2)}}^2+\|\na u\|_{L^{2r/(r-2)}}^2 \right)\\   & \le C\norm[L^r]{u}^2 \left( \norm[L^2]{\n^{1/2} E}^{2(r-3)/r} \norm[L^6]{ E}^{6/r}+  \norm[L^2]{\nabla u}^{2(r-3)/r} \norm[L^6]{\nabla  u}^{6/r} \right)\\  & \le  C (\eta) (1+\norm[L^r]{u}^{ {s} })\int\left(\n E^2+|\na u|^2 \right)dx+ \eta  \|\na E\|_{L^2}^2+ \eta \|\na u\|_{L^6}^2  ,\ea\ee where in the last inequality we  have used     \eqref{k0.2}.
It follows from  \eqref{dw},  \eqref{3.n3}, \eqref{3.n2},   and \eqref{k0.2} that \be \la{3.14-3}\ba\|\na u\|_{L^6}&\le C\|\n E\|_{L^6}+C\|\n\dot u\|_{L^2}+ C\|\na H  \|_{L^2}^{1/2}\|\na H  \|_{H^1}^{1/2}\\&\le C\|\n^{1/2} E\|_{L^2}+C\|\na E\|_{L^2}+C\|\n\dot u\|_{L^2}+C \|\na H  \|_{L^2}^{1/2}\|\na H  \|_{H^1}^{1/2}.\ea\ee
Putting this and \eqref{k0.9} into \eqref{p329}, and  choosing $\eta$ suitably small, we obtain after using Gronwall's inequality,   \eqref{p23},
 \eqref{3.3-1},  and  \eqref{1.10-1-1}   that
 \be\la{o3.12}\ba &\sup_{0\leq t\leq T}
\int\left(\r E^2 +|\nabla u|^2+|\nabla H|^2\right)
dx\\&+\int_{0}^{T}\int\left(|\nabla E|^2+\r|\dot{u}|^2+|u|^2|\na u|^2+|H_t|^2+|\na^2H|^2\right)dxdt \leq C.\ea
\ee

Finally,  \eqref{a1}$_1$ implies that \be\la{cc.2} (\n-\ti\n)_t+\div((\n-\ti\n)u)+ \ti \n\div u=0.\ee
   Multiplying \eqref{cc.2} by $\n-\ti\n$ and integrating the resulting equation over $\om ,$ we obtain after using \eqref{1.10-1-1} that \bnn (\|\n-\ti\n\|_{L^2}^2)'(t)\le C\|\n-\ti\n\|_{L^2}^2+C\|\na u\|_{L^2}^2,\enn which together with  \eqref{o3.12} and the following simple fact that \bnn \|\na \t\|_{L^2} \le C\|\na E\|_{L^2}+  C\||u||\na u|\|_{L^2},\enn
directly  gives \eqref{3.4}.
    The proof of
Lemma \ref{lem3.2} is  completed.

Finally,  the following Lemma  \ref{lem3.7} will deal with the
higher order estimates of the solutions which are needed to
guarantee the extension of local strong solution to be a global one
under the conditions  \eqref{1.7}-\eqref{1.9} and \eqref{1.10-1-1}.
\begin{lemma}\label{lem3.7}
Under the condition \eqref{1.10-1-1},  it holds that  for $0\leq
T<T^*$,
\begin{equation}\label{3.54}
\sup_{0\leq t\leq T}(\|\r-\ti\r\|_{H^1\cap
W^{1, \ti q}}+\|\nabla{u}\|_{H^1}+ \|\na\t\|_{H^1}+ \|H\|_{H^2})\leq C.
\end{equation}
\end{lemma}

{\it Proof.}
First, it  follows from \eqref{3.14-3}, \eqref{3.4}, and  \eqref{3.3-1}   that
\bnn  \ba \|\na u\|_{L^6}&\le C+C\|\n\dot u\|_{L^2} +C\|\na\te\|_{L^2}+C\||u||\na u|\|_{L^2}+C  \|\na H  \|_{H^1}^{1/2}\\ &\le C+C\|\n\dot u\|_{L^2} +C\|\na\te\|_{L^2}+C\|u\|_{L^6}\|\na u \|_{L^2}^{1/2}\|\na u \|_{L^6}^{1/2}+C \|\na H  \|_{H^1}^{1/2}\\ &\le C+C\|\n\dot u\|_{L^2} +C\|\na\te\|_{L^2} +\f12\|\na u \|_{L^6} +C\|\na H  \|_{H^1}^{1/2},\ea\enn which implies
\be\la{b3.39} \ba \|\na u\|_{L^6}&\le C+C\|\n\dot u\|_{L^2} +C\|\na\te\|_{L^2} +C \|\na H\|_{H^1}^{1/2}.\ea\ee

Then, it follows from the standard $L^2$-estimate of $(\ref{a1})_4$, \eqref{3.4},  and  \eqref{3.3-1}   that\bnn\ba \|\na^2 H\|_{L^2}&\le C\|H_t\|_{L^2}+C\||u||\na H|\|_{L^2}+C\||\na u||H|\|_{L^2}\\&\le C\|H_t\|_{L^2}+C\|u\|_{L^6}\|\na H\|^{1/2}_{L^2}\|\na^2 H\|^{1/2}_{L^2} + C\|H\|_{L^4}\|\na u\|_{L^4} \\&\le C\|H_t\|_{L^2}+\f12 \|\na^2 H\|_{L^2}+C  \|\na u\|_{L^4}+C,\ea\enn
which together with \eqref{3.4}   and \eqref{3.3-1} implies \be\la{a3.39}\ba \|  H\|_{H^2} \le   C\|H_t\|_{L^2} +C  \|\na u\|_{L^4}+C.\ea\ee
  Holder's inequality, along with  \eqref{b3.39}  and \eqref{3.4}, gives \be\la{u4} \|\na u\|_{L^4} \le \|\na u\|_{L^2}^{1/4}\|\na u\|_{L^6}^{3/4}\le C+C\|\n\dot u\|_{L^2}^{3/4}+C\|\na\te\|^{3/4}_{L^2} +C \|\na H\|_{H^1}^{3/8},\ee
which combined with  \eqref{a3.39} and   \eqref{b3.39} shows  \be \la{c3.39} \|\na u\|_{L^6}+\|  H\|_{H^2}\le C\|\n\dot u\|_{L^2} + C\|H_t\|_{L^2}+C\|\na\te\|_{L^2} +C.\ee

Then, similar to \eqref{k0.2}, we have \be\la{po0}\|\te\|_{L^6}\le C\|\n^{1/2} \te\|_{L^2}+C\|\na \te\|_{L^2}\le C+C\|\na \te\|_{L^2},\ee
which together with the
standard $L^2$-estimate of $(\ref{a1})_3$ and \eqref{3.4} gives
\bnn \ba &
\|\nabla \t\|_{H^1}^2\\ &\leq C+C\|\na\te\|_{L^2}^2+C\int\r\dot{\t}^2dx +C\int\r^2\t^2|\nabla
u|^2dx+C\|\nabla u\|_{L^4}^4  + C\|\nabla H\|_{L^4}^4  \\ &\leq C+C\|\na\te\|_{L^2}^2+C\int\r\dot{\t}^2dx +C\|\na u\|_{L^2}^2 \| \te\|^2_{L^\infty} +C\|\nabla u\|_{L^4}^4 +C \|\nabla H\|_{H^1}^4\\ &\leq C +C\int\r\dot{\t}^2dx+C\|\na \te\|_{L^2}^2+\f12 \| \na \te\|^2_{H^1} +C\|\nabla u\|_{L^4}^4 +C \|\nabla H\|_{H^1}^4 .\ea
\enn Combining this with  \eqref{a3.39}  shows
\be\label{3.44}
\|\nabla \t\|_{H^1}^2   \leq C\int\r\dot{\t}^2dx  +C\|\na\te\|^2_{L^2}+ C\|\nabla u\|_{L^4}^4 +C\| H_t\|_{L^2}^4 +C.
\ee

Next,  we claim that we have the
  following estimates on both $\dot u$ and $\dot \te,$   \eqref{3.27} and \eqref{3.46}, whose proofs are similar to those in \cite{hlx4,z} and can be found in   \ref{app}:
 \begin{equation}\label{3.27}
\sup_{0\leq t\leq
T}\int\left(|\nabla\t|^2+\r|\dot{u}|^2+|H_t|^2\right)dx+\int_{0}^{T}\int\left(\r\dot{\t}^2+|\nabla\dot{u}|^2
+|\na H_t|^2\right)dxdt\leq
C,
\end{equation}
\begin{equation}\label{3.46}
\sup_{0\leq t\leq
T}\|\r^{1/2}\dot{\t}\|_{L^2}^2 +\int_{0}^{T}\|\nabla{\dot{\t}}\|_{L^2}^2dt\leq
C.
\end{equation}

Then, the combination of    \eqref{c3.39}--\eqref{3.46}  with  \eqref{u4} leads to
\begin{equation}\label{3.55}
\sup_{0\leq t\leq
T} \left(\|\na u\|_{L^6} +\|  H \|_{H^2}+\|\t\|_{L^6}+\|\na\t\|_{H^1} \right)\leq C.
\end{equation} For $2\leq p\leq \tilde q,$  direct calculations show that
 \be \ba\label{a3.57}
\f{d}{dt}\|\nabla\r\|_{L^p} \leq& C(1+\|\nabla
u\|_{L^\infty})\|\nabla\r\|_{L^p}+C\|\nabla^2 u\|_{L^p}.
\ea\ee
For  $v \in D^1_0\cap D^{2,2}\cap D^{2,\ti q}$ satisfying \eqref{lame},  it   follows from Lemma \ref{lem2.3} and \eqref{3.55} that
\be\la{3.42}\ba \|\na v\|_{L^\infty}&\le C\left( 1+\log \left(e+\|\na (\n\te)\|_{L^{\ti q}}\right)\|\n\te\|_{L^\infty}+\|\n\te\|_{L^2}\right)\\ &\le C \log \left(e+\|\na \n\|_{L^{\ti q}}  \right).\ea\ee
Then, for $w\triangleq u-v \in D^1_0\cap D^{2,2}\cap D^{2,\ti q}$ satisfying \eqref{wl1}, applying the standard $L^p$-estimate to  \eqref{wl1}, along with \eqref{3.55}, gives
\bnn\ba \|\na^2w\|_{L^6}&\le C\|\n\dot u\|_{L^6} +C\||H||\na H| \|_{L^6}\\ &\le  C\|\na\dot u\|_{L^2}+C ,\ea\enn
which together with  \eqref{3.n2},    \eqref{3.55}, and \eqref{3.27} shows \bnn \|\na w\|_{L^\infty}\le C+ C\|\na\dot u\|_{L^2}. \enn The combination of this with \eqref{3.42}  gives \be\la{3.43} \|\na u\|_{L^\infty}\le C \log \left(e+\|\na \n\|_{L^{\ti q}}  \right)+ C\|\na\dot u\|_{L^2}.\ee
Applying the standard $L^p$-estimate to \eqref{a1}$_2$ leads to \be\ba
 \la{3.59}\|\na^2 u\|_{L^p}  &\le   C\left(\|\n\dot u\|_{L^p}+
  \||H||\na H|\|_{L^p}+
  \|\nabla P\|_{L^p}\right)\\ & \le   C\left(\|\n\dot u\|_{L^p}+    \|\nabla\n\|_{L^p}\right)+C\\ & \le   C\left(1+\|\na\dot u\|_{L^2}+\|\nabla\n\|_{L^p}\right) ,\ea\ee
due to   \eqref{3.27} and \eqref{3.55}.
 Substituting \eqref{3.59} and \eqref{3.43}  into \eqref{a3.57} yields that
\begin{equation}\la{a3.60}
f'(t)\leq  Cg(t)f(t)\ln{f(t)},
\end{equation} where
\bnn
 f(t)\triangleq e+\|\nabla\r\|_{L^{\ti q}}, \,\, g(t)\triangleq
1+\|\nabla\dot{u}\|_{L^2} .
\enn
It thus follows from  \eqref{a3.60}, \eqref{3.27},  and    Gronwall's inequality that
 \begin{equation}\label{3.64}
\sup_{0\leq t\leq T}\|\nabla\r\|_{L^{\ti q}}\leq C,
\end{equation}
which, along with \eqref{3.43} and \eqref{3.27},
 directly gives
\begin{equation}\label{3.65}
\int_{0}^{T}\|\nabla u\|^2_{L^\infty}dt\leq C.
\end{equation}
Taking $p=2$ in \eqref{a3.57}, we get by using \eqref{3.65}, \eqref{3.59}, \eqref{3.27}, and Gronwall's inequality  that
\begin{equation}\label{3.66}
\sup_{0\leq t\leq T}\|\nabla\r\|_{L^2}\leq C,
\end{equation}
which together with  \eqref{3.59},
\eqref{3.55},   and
\eqref{3.27}   yields that
\begin{eqnarray*}\label{3.67}
\sup_{0\leq t\leq T}\|\nabla^2u\|_{L^2}\leq C\sup_{0\leq t\leq
T}\left(\|\r\dot{u}\|_{L^2}+\|\nabla\r\|_{L^2}
+\|\nabla\t\|_{L^2} + \|H\cdot\na H\|_{L^2}\right) \leq C.
\end{eqnarray*}
This combined  with \eqref{3.64},  \eqref{3.66}, \eqref{3.55}, and
\eqref{3.4} finishes the proof of Lemma \ref{lem3.7}.

  Now we are in a position to prove Theorem  \ref{thm1.1}.

{\it Proof  of Theorem \ref{thm1.1}.}
Suppose that \eqref{1.10} were false, that is, \eqref{1.10-1-1} holds.
 Note that the generic constant  $C$ in Lemma
 \ref{lem3.7} remains uniformly bounded for all
$T<T^\ast$, so the functions
$(\r,u,\t, H)(x,T^\ast)\triangleq\lim\limits_{t\rightarrow
T^\ast}(\r,u,\t, H)(x,t)$ satisfy the conditions imposed on the initial
data \eqref{1.7} at the time $t=T^\ast$. Furthermore, standard
arguments yield that $\r\dot{u},\r\dot{\t}\in C([0,T];L^2)$, which
implies \bnn (\r\dot{u},\r\dot{\t})(x,T^\ast)=\lim_{t\rightarrow
T^\ast}(\r\dot{u},\r\dot{\t})\in L^2. \enn Hence, \bnn
&&-\m\Delta{u}-(\m+\l)\nabla\mbox{div}u+R\nabla(\r\t)
-(\mbox{curl $H$})\times H|_{t=T^\ast}=\sqrt{\r}(x,T^\ast)g_1(x),\nonumber\\
&&\k\Delta\t+\f\m2|\nabla{u}+(\nabla{u})^{tr}|^2+\l(\mbox{div}u)^2
+\nu|\mbox{curl $H$}|^2|_{t=T^\ast}=\sqrt{\r}(x,T^\ast)g_2(x)\nonumber,
\enn with \bnn  g_1(x)\triangleq
\begin{cases}
\r^{-1/2}(x,T^\ast)(\r\dot{u})(x,T^\ast),&
\mbox{for}~~x\in\{x|\r(x,T^\ast)>0\},\\
0,&\mbox{for}~~x\in\{x|\r(x,T^\ast)=0\},
\end{cases}
\enn and \bnn g_2(x)\triangleq
\begin{cases}
\r^{-1/2}(x,T^\ast)(c_v\r\dot{\t}+R\r\t\mbox{div}u)(x,T^\ast),&
\mbox{for}~~x\in\{x|\r(x,T^\ast)>0\},\\
0,&\mbox{for}~~x\in\{x|\r(x,T^\ast)=0\},
\end{cases}
\enn satisfying $g_1,g_2\in L^2$ due to \eqref{3.27}, \eqref{3.46},
and \eqref{3.54}. Thus, $(\r,u, \t,H)(x,T^\ast)$ also satisfies \eqref{1.8} and
\eqref{1.9}. Therefore, one can take $(\r,u, \t,H)(x,T^\ast)$ as
the initial data and apply Lemma \ref{lem2.1} to extend the local
strong solution beyond $T^\ast$. This contradicts the assumption on
$T^{\ast}$. We thus finish the proof of Theorem \ref{thm1.1}.

\appendix
\renewcommand\thesection{\appendixname~\Alph{section}}
\renewcommand\theequation{\Alph{section}.\arabic{equation}}

\section{\la{app} Proofs of  \eqref{3.27} and \eqref{3.46}.}
 The proofs of  \eqref{3.27} and \eqref{3.46} are a direct combination of those of Lemma 4.1 and (4.28) in \cite{hlx4} with that of (3.24) in \cite{z}. We sketch them here for completeness.

 First, it follows from \eqref{3.4} and  \eqref{c3.39} that
\be\label{3.4a}\ba &\sup_{0\leq t\leq T}
\int \r \te^2
dx +\int_{0}^{T} \left(\|\nabla \te\|_{L^2}^2+\|\nabla u\|_{L^6}^2\right) dt \leq C.\ea
\ee
 Applying $\dot{u}_j[\partial_t+\mbox{div}(u\cdot)]$ to $ (\ref{a1})_2^j$ and integrating the resulting
equality over $\om $ give
\be\label{3.23}
\ba
\f12\f{d}{dt}\int\r|\dot{u}|^2dx =&-\int\dot{u}_j
[\partial_jP_t+\mbox{div}(u\partial_jP)]
dx+\m\int\dot{u}_j[\pa_t\Delta{u}_j +\mbox{div}(u\Delta{u}_j)]dx \\
&+(\m+\l)\int\dot{u}_j[\partial_j\mbox{div}u_t+\mbox{div}(u\partial_j{\rm div}u)]dx \\
& -\frac{1}{2}\int\dot{u}_j[\p_t\p_j|H|^2 + \mbox{div}(u\p_j|H|^2)]dx\\
& +
\int\dot{u}_j[\p_t\p_i(H^i H_j)+\mbox{div}(u\p_i(H^i H_j))]dx\\
=&\sum_{i=1}^{5}N_i.
\ea
\ee
We get  after integration by parts and
using  the equation $(\ref{a1})_1$ that
\be\la{3.24} \ba
N_1 & = - \int   \dot{u}_j[\partial_jP_t + \div (\p_jPu)]dx \\
& =R\int   \p_j\dot{u}_j\left(\n \dot\te-\r u\cdot\na\te-\te
 u\cdot\na\n-\te  \n\div u\right) dx + \int  \p_k\dot{u}_j\p_jPu_k dx  \\
& =R\int   \p_j\dot{u}_j\left(\n \dot\te
 -\te  \n\div u\right)dx
 +\int  P \div\dot{u}  \div u dx
   -  \int    P\pa_k\dot{u}_j\p_ju_k dx
 \\&\le \frac{ \mu}{8} \|\nabla\dot{u}\|_{L^2}^2
 +C    \|\n \dot\te\|_{L^2}^2
+C\int\n^2\te^2|\na u|^2dx   \\&\le \frac{ \mu}{8} \|\nabla\dot{u}\|_{L^2}^2
 +C    \|\n \dot\te\|_{L^2}^2
+C\|\n\t\|_{L^2}^{1/2}\|\te\|_{L^6}^{3/2}\|\na u\|_{L^4}^2 \\&\le \frac{ \mu}{8} \|\nabla\dot{u}\|_{L^2}^2
 +C    \|\n \dot\te\|_{L^2}^2
+C \|\na\te\|_{L^2}^{4}+C\|\na u\|_{L^4}^4+C,\ea \ee where in the last inequality  we have used \eqref{po0}.
Integration by parts leads to
\be\label{3.25} \ba N_2 & =  \mu\int
 \dot{u}_j[\pa_t\triangle u_j
 + \div (u\triangle u_j)]dx \\
& = - \mu\int  \left(\p_i\dot{u}_j(\p_iu_j)_t +
\triangle u_ju\cdot\nabla\dot{u}_j\right)dx \\
& = -  \mu\int \left(|\nabla\dot{u}|^2 -
\p_i\dot{u}_ju_k\p_k\p_iu_j - \p_i\dot{u}_j\p_iu_k\p_ku_j +
\triangle u_ju\cdot\nabla\dot{u}_j\right)dx \\
& = - \mu\int  \left(|\nabla\dot{u}|^2 + \p_i\dot{u}_j
\p_iu_j\div u - \p_i\dot{u}_j\p_iu_k\p_ku_j - \p_iu_j\p_iu_k\p_k\dot{u}_j
\right)dx \\
&\le -\frac{7\mu}{8} \int |\nabla\dot{u}|^2dx  + C \int
 |\nabla u|^4dx  . \ea \ee
Similarly, we have
\begin{eqnarray}\label{3.26}
N_3&\leq&
-\f78(\m+\l)\|\mbox{div}\dot{u}\|_{L^2}^2+C\int|\nabla{u}|^4dx.
\end{eqnarray}
Integration by parts together with \eqref{3.4} and \eqref{3.3-1} shows
\be\la{n-4}
\ba
|N_4|
& \le C\norm[L^2]{\nabla\dot{u} }( \||H||H_t|\|_{L^2}+\||u||H||\na H|\|_{L^2})\\& \le C\norm[L^2]{\nabla\dot{u} }( \|H\|_{L^6}\|H_t\|_{L^2}^{1/2}\|H_t\|_{L^6}^{1/2} +\|u\|_{L^6}\|H\|_{L^6}\|\na H\|_{L^6})\\
& \le \ep\norm[L^2]{\nabla\dot{u}}^2+ \eta \norm[L^2]{\na H_t}^2+ C(\ve,\eta) \norm[L^2]{ H_t}^2 +
C(\ve)\norm[L^2]{\nabla^2 H}^2 .
\ea
\ee
Similarly, we also have
\be\la{n-5}
\ba|N_5| \le\ep\norm[L^2]{\nabla\dot{u}}^2+ \eta \norm[L^2]{\na H_t}^2+ C(\ve,\eta) \norm[L^2]{H_t}^2+
C(\ve)\norm[L^2]{\nabla^2 H}^2 .
\ea
\ee
Substituting \eqref{3.24}-\eqref{n-5} into \eqref{3.23}, we obtain after choosing $\ve$ suitably small that
\begin{equation}\label{3.22}\ba
  \f{d}{dt}\int\r|\dot{u}|^2dx+\m\|\nabla\dot{u}\|_{L^2}^2  &\leq
C\int\r\dot{\t}^2dx+C\eta\norm[L^2]{\nabla H_t}^2+C( \eta) \norm[L^2]{H_t}^2\\&\quad +
C \norm[L^2]{\na^2 H}^2 +C\|\nabla\t\|_{L^2}^4+C\|\nabla{u}\|_{L^4}^4+C.\ea
\end{equation}

Next, multiplying $(\ref{a1})_3$ by $  \dot\te$ and integrating
the resulting equality over $\om $ yield  that
 \be\la{e1} \ba   &\frac{\ka
{ }}{2}\left( \|\na\te\|_{L^2}^2\right)_t +c_v\int\rho|\dot{\te}|^2dx  \\
&=-\ka \int\na\te\cdot\na(u\cdot\na\te)dx
+\lambda \int  (\div u)^2\dot\te dx\\&
\quad+2\mu \int |\mathfrak{D}(u)|^2\dot\te dx -R \int\n\te \div
u\dot\te dx +\nu\int|\mbox{curl}H|^2\dot{\te} dx \\&\triangleq \sum_{i=1}^5I_i . \ea\ee

We estimate each $I_i (i=1,\cdots,5)$ as follows:

First, it follows from \eqref{3.44} and  (\ref{3.4}) that
\be\la{e2} \ba  |I_1|&\le C  \int|\na u||\na\te|^2dx\\
&\le C  \|\na u\|_{L^2}\|\na\te\|^{1/2}_{L^2}
\|\na\te\|^{3/2}_{L^6} \\
&\le \de  \|\na^2\te\|^2_{L^2} +C(\de)
\|\na\te\|^2_{L^2}  \\
&\leq C\de\int\r\dot{\t}^2dx +C(\de)\|\na\te\|^2_{L^2}+  C\|\nabla u\|_{L^4}^4 +C\| H_t\|_{L^2}^4 +C .\ea\ee

Next, integration by parts yields that, for any $\eta \in (0,1],$
\be\la{e3}\ba I_2 =&\lambda \int (\div u)^2 \te_t
dx+\lambda \int (\div u)^2u\cdot\na\te
dx\\=&\lambda \int\left( (\div u)^2 \te
\right)_t dx-2\lambda  \int \te \div u \div (\dot u-u\cdot\na u)
dx +\lambda \int (\div u)^2u\cdot\na\te
dx \\=&\lambda
\left(\int (\div u)^2 \te dx\right)_t-2\lambda \int \te \div u
\div \dot udx\\&+2\lambda \int \te \div u \pa_i u_j\pa_j  u_i dx
+ \lambda \int u \cdot\na\left(\te   (\div u)^2 \right)dx
  \\ \le &\lambda\left( \int (\div u)^2 \te dx\right)_t  +C\| \te\|_{L^6} \|\na u\|_{L^2}^{1/3}\|\na u\|_{L^4}^{2/3} \left(\|\na \dot u\|_{L^2}+\|\na u\|_{L^4}^2\right)\\ \le &\lambda\left( \int (\div u)^2 \te dx\right)_t +\eta\|\na \dot u\|_{L^2}^2  +C(\eta)\|\na u\|_{L^4}^4 +C\|\na\te\|_{L^2}^4+C,
 \ea\ee  where in the last inequality  we have used \eqref{po0}.

 Then, similar to (\ref{e3}), we have that, for
any $\eta  \in (0,1],$\be \la{e5}\ba I_3\le& 2\mu\left( \int
|\mathfrak{D}(u)|^2 \te dx\right)_t +\eta\|\na \dot u\|_{L^2}^2  +C(\eta)\|\na u\|_{L^4}^4 +C\|\na\te\|_{L^2}^4+C. \ea\ee

Next, it follows from \eqref{3.4} and \eqref{po0}    that
 \be\la{e39}\ba
 |I_4|   & \le   C   \|\n^{1/2}  \dot\te\|_{L^2} \|\n^{1/2} \te\|_{L^2}^{1/4}  \| \te\|_{L^6}^{3/4}\|\na u\|_{L^4} \\ & \le \de\int\r\dot{\t}^2dx +C(\de)\|\na\te\|^4_{L^2} +C \|\na u\|_{L^4}^4  +C(\de),  \ea\ee
and  that
\be\la{i-5}
\ba
I_5 & = \nu\int|\mbox{curl} H|^2\te_tdx + \nu\int|\mbox{curl} H|^2u\cdot\nabla\te dx\\
& =\nu\frac{d}{dt}\int|\mbox{curl} H|^2\te dx - 2\nu\int\te\mbox{curl} H\cdot\mbox{curl} H_tdx + \nu\int|\mbox{curl} H|^2u\cdot\nabla\te dx\\
& \le\nu\frac{d}{dt}\int|\mbox{curl} H|^2\te dx +
C\norm[L^6]{\te}\|\na H\|_{L^2}^{1/2}\|\na H\|_{L^6}^{1/2}\norm[L^2]{\nabla H_t} \\& \quad+C\|\na H\|_{L^6}^2\|u \|_{L^6} \norm[L^2]{\nabla \te}  \\
& \le \nu\frac{d}{dt}\int|\mbox{curl} H|^2\te dx +
\eta\norm[L^2]{\nabla H_t}^2 + C(\eta)\left(1+\norm[L^2]{\nabla \te}^2\right)\left(1+\|\na^2H\|_{L^2}^2\right) .
\ea
\ee

Substituting \eqref{e2}-\eqref{i-5} into \eqref{e1}, we obtain after choosing $ \de$ suitably small that, for any $\eta\in (0,1],$
 \be\la{ee1} \ba   &\f{d}{dt}\int\Psi dx+c_v\int\rho|\dot{\te}|^2dx  \\
&\le C(\eta )\left(1+\norm[L^2]{\nabla \te}^2\right)\left(1+\|\na^2H\|_{L^2}^2+\norm[L^2]{\nabla \te}^2\right)+C\eta\|\na \dot u\|_{L^2}^2\\&\quad+C\eta\norm[L^2]{\nabla H_t}^2+  C\|\nabla u\|_{L^4}^4   +C\|H_t\|_{L^2}^4+C   , \ea\ee
where \be \la{po}\Psi\triangleq\k |\nabla\t|^2-2\t\left[\l(\mbox{div}
u)^2+2\m|\mathfrak
{D}(u)|^2+ {\nu}|\mbox{curl} H|^2\right].\ee

Next,  differentiating $(\ref{a1})_4$ with respect to $t$ and multiplying the resulting equations by $H_t,$ we obtain after   integration   by parts and using \eqref{3.3-1} and \eqref{3.4}  that
\bnn
\ba
& \f12\frac{d}{dt}\int|H_t|^2dx +\nu\int|\nabla H_t|^2dx \\
&  \le C\left(\||u_t|| H|\|_{L^2}+\| |u||H_t|\|_{L^2}\right)\|\na H_t\|_{L^2}\\&\le C\left(\||\dot u|| H|\|_{L^2}+\||u||\na u|| H|\|_{L^2}+\| u\|_{L^6}\|H_t\|_{L^2}^{1/2}\|  H_t\|_{L^6}^{1/2}\right)\|\na H_t\|_{L^2}\\&\le C\left(\|\dot u\|_{L^6}\| H\|_{L^3}+\|u\|_{L^6}\|\na u\|_{L^4}\| H\|_{L^{12}}+ \|H_t\|_{L^2}^{1/2}\|\na H_t\|_{L^2}^{1/2}\right)\|\na H_t\|_{L^2}\\&\le  \frac{\nu}{2}\|\na H_t\|_{L^2}^2 +C\|\na\dot u\|_{L^2}^2+ C\|\na u\|_{L^4}^2+C\|H_t\|_{L^2}^2 ,
\ea
\enn
 which implies \be\la{ht-1} \frac{d}{dt}\int |H_t|^2dx +\nu \int|\nabla H_t|^2dx\le C\|\na\dot u\|_{L^2}^2+ C\|\na u\|_{L^4}^2+C\|H_t\|_{L^2}^2 .\ee

Finally, adding \eqref{3.22} multiplied by $ \eta^{1/4}$ and  \eqref{ht-1} by $ \eta^{1/2}$  to  \eqref{ee1}, we obtain after choosing $\eta$ suitably small and using \eqref{u4} that
 \be\la{3.46a}\ba &2\frac{d}{dt}\int (\Psi+\eta^{1/2}|H_t|^2+\eta^{1/4}\n|\dot u|^2)dx\\&\quad +\int \left(c_v\n|\dot \te|^2+ \nu{\eta^{1/2}}|\na H_t|^2+\mu\eta^{1/4}|\na\dot u|^2\right)dx \\&\le   C(\eta )\left(1+\norm[L^2]{\nabla \te}^2\right)\left(1+\|\na^2H\|_{L^2}^2+\norm[L^2]{\nabla \te}^2\right)\\&\quad+  C(\eta )\|\r^{1/2}\dot{u}\|^{4}_{L^2} +C(\eta)\|H_t\|_{L^2}^4 . \ea\ee

Noticing that  \eqref{po}, \eqref{3.4},  \eqref{c3.39},  and  \eqref{po0} lead to  \bnn\ba &2\int (\Psi+\eta^{1/2}|H_t|^2+\eta^{1/4}\n|\dot u|^2)dx\\&\ge 2\k \|\na\te\|_{L^2}^2-C\|\t\|_{L^6}\|\nabla{u}\|^{3/2}_{L^2}\|\nabla{u}
\|^{1/2}_{L^6}-
C\|\t\|_{L^6}\|\nabla H\|_{L^2}^{3/2}\|\na H\|_{H^1}^{1/2}\\&\quad +2\int ( \eta^{1/2}|H_t|^2+\eta^{1/4}\n|\dot u|^2)dx\\&\ge   \k \|\na\te\|_{L^2}^2    + \int ( \eta^{1/2}|H_t|^2+\eta^{1/4}\n|\dot u|^2)dx-C (\eta),\ea\enn
 we directly obtain \eqref{3.27} after   using Gronwall's inequality,   \eqref{1.8}, \eqref{3.46a},   \eqref{3.4},  and  \eqref{3.4a}.

   It remains to prove \eqref{3.46}. First,  it follows from \eqref{c3.39}--\eqref{3.27} that
\be\la{hhh}
\sup_{0\le t\le T}(\|\te\|_{L^6}+\|\na u\|_{L^2\cap L^6}+\| H\|_{H^2})
+\int_0^T \|\na^2\t\|_{L^2}^2 dt\le C.
\ee
Similar to \eqref{k0.2}, we have \be\la{k03} \|\dot \te\|_{L^6}\le C\|\n^{1/2}\dot \te\|_{L^2}+C\|\na\dot \te\|_{L^2}.\ee

Next, applying the operator $\pa_t+\div(u\cdot) $ to (\ref{a1})$_3$ leads to \be\la{3.96}\ba
&c_v\n \left(\pa_t\dot \te+u\cdot\na\dot \te\right)\\
&=\ka \Delta \dot \te+\ka\left( \div u\Delta \te  -
\pa_i\left(\pa_iu\cdot\na \te \right)- \pa_iu\cdot\na \pa_i\te
\right)\\&\quad +\left( \lambda (\div u)^2+2\mu |\mathfrak{D}(u)|^2\right)\div u
+R\n \te  \pa_ku_l\pa_lu_k
\\&\quad -R\n \dot\te \div u-R\n \te\div \dot u
 +2\lambda \left( \div\dot u-\pa_ku_l\pa_lu_k\right)\div u\\&\quad
+ \mu (\pa_iu_j+\pa_ju_i)\left( \pa_i\dot u_j+\pa_j\dot
u_i-\pa_iu_k\pa_ku_j -\pa_ju_k\pa_ku_i\right)\\
&\quad + \nu\left[\p_t|\mbox{curl}H|^2 + \mbox{div}(|\mbox{curl}H|^2u)\right].\ea\ee
Multiplying
(\ref{3.96}) by $\dot \te,$  we obtain after integration by parts
and using \eqref{hhh}, \eqref{3.27}, and \eqref{k03} that \bnn\la{3.99}\ba &
\frac{c_v}{2}\left(\int \n |\dot\te|^2dx\right)_t + \ka
\|\na\dot\te\|_{L^2}^2 \\&\le
 C  \int|\na u|\left(|\na^2\te||\dot\te|+ |\na \te| |\na
\dot\te|\right)dx +C   \int|\na
u|^2|\dot\te|\left(|\na u|+\te \right)dx
 \\&\quad +C   \int  \n  |\dot
\te|^2|\na u| dx   +C \int  \n \te  |\na\dot u| |\dot
\te|dx +C   \int  |\na u| |\na\dot u| |\dot
\te|dx\\
&\quad+ C\int\left(|\na H||\na H_t||\dot{\t}|+|\na H|^2|u||\na\dot{\t}|\right)dx
\\&\le C  \|\na u\|_{L^3}\|\na\te\|_{H^1}\left(\|\dot\te\|_{L^6}+\|\na\dot \te\|_{L^2} \right)+C  \|\na u\|_{L^3}^2
\| \dot\te\|_{L^6} \left(\|\na
u\|_{L^6}+\| \te\|_{L^6}\right) \\
&\quad+C  \|\na u\|_{L^3} \|\n\dot\te\|_{L^2}
\|\dot\te\|_{L^6}
      +C \|\n^{1/2}\te\|^{1/2}_{L^2}\|\te\|_{L^6}^{1/2}
 \|\na\dot u\|_{L^2} \|\dot\te\|_{L^6}  \\&\quad+C    \|\na
u\|_{L^3} \|\na\dot u\|_{L^2} \|\dot\te\|_{L^6} + C\|\na H\|_{L^3}\|\na H_t\|_{L^2}\|\dot{\t}\|_{L^6} \\&\quad+
C\|\na H\|_{L^6}^2\|u\|_{L^6}\|\na\dot{\t}\|_{L^2}
\\ & \le\frac{\ka}{2}\|\na\dot\te\|_{L^2}^2+C\|\na^2\t\|_{L^2}^2
 +C\|\n^{1/2}\dot\t\|_{L^2}^2 +C\|\na\dot u\|_{L^2}^2 +
 C\|\na H_t\|_{L^2}^2+C,   \ea\enn
 which combined with Gronwall's inequality, \eqref{1.9},  \eqref{hhh}, and \eqref{3.27} directly gives  \eqref{3.46}.

\end{document}